\documentclass[11pt]{article}

\usepackage{amsmath}
\usepackage{graphicx}
\usepackage{amssymb}
\usepackage{multicol}
\usepackage[usenames]{color}
\usepackage[sc,small]{caption}

\graphicspath{{./figs/}, {./oldfigs/}}

\definecolor{RED}{rgb}{1,0,0}
\definecolor{BLUE}{rgb}{0,0,1}
\definecolor{PURPLE}{rgb}{1,0,1}

\newcommand{\order}{\ensuremath{\mathcal{O}}}
\newcommand{\Lyap}{\ensuremath{\mathcal{F}}}
\newcommand{\sech}{\ensuremath{\mathrm{sech}}}
\newcommand{\cc}{\ensuremath{\mathrm{c.c.}}}
\newcommand{\brLzero}{\ensuremath{\mathrm{L0}}}
\newcommand{\brLone}{\ensuremath{\mathrm{L1}}}
\newcommand{\brP}{\ensuremath{\mathrm{P}}}
\newcommand{\segonepeak}{\ensuremath{\brLzero_\mathrm{1pk}}}

\begin{document}
\author{J.~Burke${}^{\dag}$\footnote{Corresponding author: { jb@math.bu.edu}}  \ and J. H. P. Dawes${}^{\ddag}$ \\
${}^{\dag}$\small Department of Mathematics and Statistics, Boston University, Boston, MA 02215, USA \\ 
${}^{\ddag}$\small Dept. of Mathematical Sciences, University of Bath, Claverton Down, Bath BA2 7AY, UK
}
\title{Localised states in an extended Swift--Hohenberg equation}
\date{\today}
\maketitle


\begin{abstract}

Recent work on the behaviour of localised states in pattern forming partial differential equations has focused on the traditional model Swift-Hohenberg equation which, as a result of its simplicity, has additional structure --- it is variational in time and conservative in space. In this paper we investigate an extended Swift-Hohenberg equation in which non-variational and non-conservative effects play a key role. Our work concentrates on aspects of this much more complicated problem. Firstly we carry out the normal form analysis of the initial pattern forming instability that leads to small-amplitude localised states. Next we examine the bifurcation structure of the large-amplitude localised states. Finally we investigate the temporal stability of one-peak localised states. Throughout, we compare the localised states in the extended Swift-Hohenberg equation with the analogous solutions to the usual Swift-Hohenberg equation. 

\end{abstract}

Keywords: homoclinic snaking, pattern formation, dissipative solitions

\section{Introduction and background} \label{sec:intro}

The study of the spontaneous emergence of patterns of activity out of homogeneous states has a long history, motivated by a wealth of examples in fluid and solid mechanics, and more recently extended to nonlinear optics, granular media, chemical reactions and mathematical biology. General surveys are given, from different viewpoints, by Cross and Hohenberg~\cite{CrossHohenberg1993}, Cross and Greenside~\cite{CrossGreenside}, Hoyle~\cite{Hoyle2006} and Pismen~\cite{Pismen2006}.

In many cases the stable patterns that appear are roughly periodic in space and extend throughout the bulk of the experimental or computational domain, growing smoothly in amplitude as a control parameter increases above a critical value. In other cases the system displays hysteresis: the pattern appears abruptly at finite amplitude as the system undergoes bifurcation, and the pattern persists as the control parameter is then reduced below the bifurcation point. This latter case is often referred to as the `subcritical' case, in contrast with the former which is the `forward' or `supercritical' case.

One of the most popular model equations for the examination of the dynamics of pattern-forming systems of this kind on a domain $\Omega \subseteq \mathbb{R}$ is the Swift-Hohenberg equation~\cite{SwiftHohenberg1977} in one spatial dimension with quadratic-cubic nonlinearity:
\begin{equation}\label{eq:intro.sh23}
 \partial_t u = r u - \left( 1 + \partial_{xx}\right)^2 u + b u^2 - u^3 \, ,
\end{equation}
where $u(x,t)$ is a scalar variable which describes the pattern forming activity, and $r$ and $b$ are real control parameters. Typical analyses~\cite{BurkeKnobloch2006} fix $b$ and treat $r$ as the primary bifurcation parameter. The trivial state $u(x,t)=0$ is linearly stable in $r<0$ and undergoes a pattern forming instability at $r=0$, where the spatial coordinate in~(\ref{eq:intro.sh23}) is scaled so that the initial instability is to perturbations with wavenumber $k=1$. In $r>0$ the trivial state is unstable to steady, spatially periodic perturbations with a range of wavenumbers surrounding $k=1$. The secondary bifurcation parameter $b$ determines the criticality  of the pattern forming instability at $r=0$: it is supercritical if $b^2 < 27/38$ and subcritical if $b^2 > 27/38$.

Equation~(\ref{eq:intro.sh23}) has several important symmetries. Firstly, it is equivariant under spatial reflections $(x,u) \rightarrow (-x,u)$. Secondly, it is equivariant under the following inversion involving the parameters: $(x,u;r,b) \rightarrow (x, -u; r, -b)$. In consequence, the behaviour of~(\ref{eq:intro.sh23}) can be fully classified by considering only the $b \geq 0$ half of the parameter plane.

We note that a common variant of~(\ref{eq:intro.sh23}) replaces the quadratic-cubic nonlinearities $b u^2 - u^3$ with the cubic-quintic combination $c u^3 - u^5$~\cite{BurkeKnobloch2007}. This results in a subcritical instability for all $c>0$, and is appropriate if there is an additional symmetry $u \rightarrow -u$ of the system, as for example in the case of Boussinesq thermal convection with identical upper and lower boundary conditions~\cite{Hoyle2006}. Of course~(\ref{eq:intro.sh23}) also includes this extra symmetry at $b=0$, and this is a popular choice to model a supercritical pattern forming system. In this article, we are interested in the more generic case in which the additional $u \rightarrow -u$ symmetry is absent.

The dynamics of~(\ref{eq:intro.sh23}) are strongly influenced by the fact that it is variational --- i.e., it can be written in the form
\begin{equation} \label{eq:intro.variational}
 \partial_t u = - \frac{ \delta \Lyap[u]}{\delta u} \, ,
\end{equation}
where the Lyapunov functional $\Lyap[ u ]$ (which we refer to as a free energy) is given by: 
\begin{equation} \label{eq:intro.Lyap}
  \Lyap[u] = \int_{\Omega} \left\{ -\frac{1}{2} r u^2 + \frac{1}{2} \left( (1 + \partial_{xx}) u \right)^2 - \frac{1}{3} b u^3 + \frac{1}{4} u^4 \right\} dx \, ,
\end{equation} 
on the bounded domain $\Omega$ with suitable boundary conditions on $\partial\Omega$ (for example, Neumann). It follows that 
\begin{equation} \label{eq:intro.dFdt}
 \frac{d}{dt} \Lyap[u] = -\int_{\Omega} (\partial_t u)^2 dx \leq 0 \, , 
\end{equation}
so that the free energy decreases in time along trajectories. Straightforward calculations show that $\Lyap[u]$ is in addition bounded from below. Hence the variational property guarantees that solutions converge to equilibria: sustained oscillations, travelling waves, and temporal chaos can not arise. The resulting equilibrium profiles $u(x)$ satisfy the time-independent version of~(\ref{eq:intro.sh23}). An equivalent description of these equilibria is as trajectories in $x$ of a fourth-order spatial dynamical system. For example, a stationary, uniform amplitude pattern of wavenumber $k$ corresponds to a periodic orbit in $x$ with period $2\pi/k$. The spatial dynamical system thus derived is conservative, and the quantity
\begin{equation}
  H = -\frac{1}{2} (r-1) u^2 + (\partial_x u)^2 - \frac{1}{2} (\partial_{xx} u)^2 + (\partial_x u) (\partial_{xxx} u) - \frac{1}{3} b u^3 + \frac{1}{4} u^4 
\end{equation}
is independent of $x$. Both of these properties --- variational in $t$ and conservative in $x$ --- aide in the analysis of~(\ref{eq:intro.sh23}).

\begin{figure}
\begin{center}
 \resizebox{6in}{!}{\includegraphics{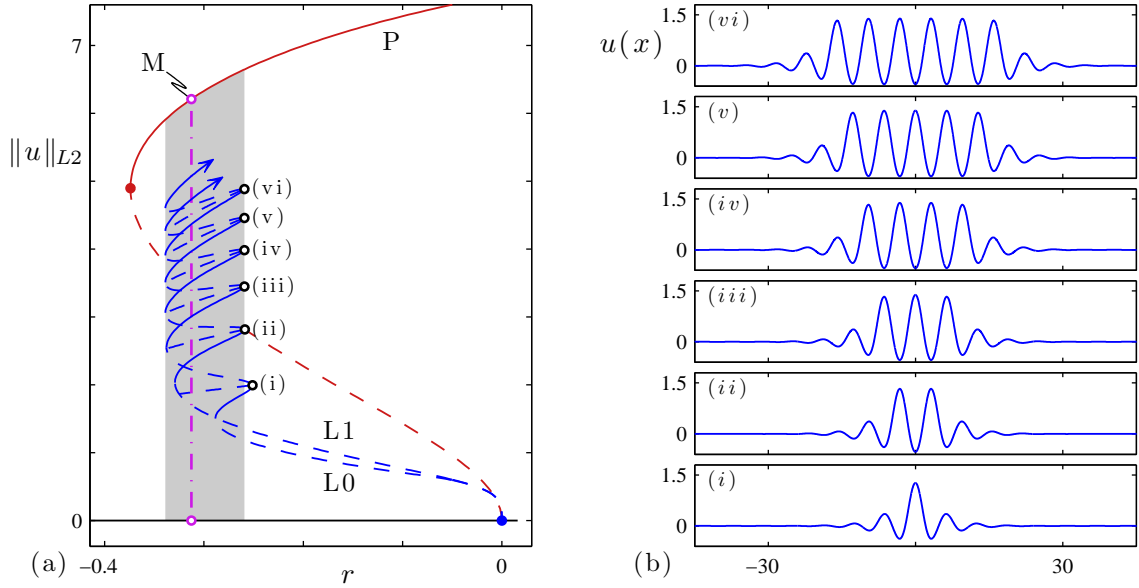}}
\end{center}
\caption{(a) Bifurcation diagram of stationary solutions to equation~(\ref{eq:intro.sh23}) at $b=1.8$, plotted in terms of the norm $||u||_{L2} = (\int_\Omega u^2(x) dx)^{1/2}$. Shading indicates the snaking region. The snaking branches $\brLzero$ and $\brLone$ include even-symmetric localised states; the arrows indicate that on $\Omega=\mathbb{R}$ the snaking continues indefinitely. The rung branches which cross-link the snaking branches are also shown. The branch $\brP$ of spatially periodic patterns satisfies $H=0$, and includes the Maxwell point $\mathrm{M}$ at which $\Lyap=0$. The norm of solutions on $\brP$ is rescaled so that this branch can be displayed on the same scale as the branches of localised states. Solid/dashed curves indicate stable/unstable solutions. (b) Profiles from several saddle-node bifurcations of the snaking branches --- profiles (i), (iii) and (v) are from $\brLzero$, and profiles (ii), (iv) and (vi) are from $\brLone$.}
\label{fig:sh23}
\end{figure}

On spatially extended domains $\Omega = \mathbb{R}$, the Swift-Hohenberg equation~(\ref{eq:intro.sh23}) exhibits a variety of equilibria beyond the uniform amplitude patterns mentioned above. In the subcritical regime, this includes stationary, spatially localised states for which the amplitude of the pattern decays to $u \rightarrow 0$ as $x \rightarrow \pm \infty$. In the spatial dynamical system, these localised profiles correspond to orbits homoclinic in $x$ to the trivial fixed point $u=0$. Recent short reviews of the generation of localised states in this context are given by Knobloch~\cite{Knobloch2008} and Dawes~\cite{Dawes2010}. The bifurcation analysis of the initial pattern-forming instability at $r=0$ (equivalently referred to in the context of fourth-order dynamical systems as a reversible 1:1 resonance or a Hamiltonian-Hopf bifurcation) was initiated by Iooss and Perou\`{e}me~\cite{IoossPeroueme1993} and developed by later authors, notably Woods and Champneys~\cite{WoodsChampneys1999}, Coullet et al.~\cite{Coullet2000}, and Burke and Knobloch~\cite{BurkeKnoblochDCDS}. The normal form analysis shows that in the subcritical regime a family of steady, small amplitude localised states of the form 
\begin{equation}
 u(x) \sim \sqrt{-r} \, \sech(x \sqrt{-r}/2) \, \cos(x+\phi) + \order(r)
\end{equation}
bifurcates from $u=0$ into $r<0$, along with the uniform amplitude patterns. Terms present in the Swift-Hohenberg equation that are manifested beyond all algebraic orders in the normal form break the $S^1$ normal form symmetry associated with $\phi$ and select $\phi=0$ and $\phi=\pi$ as the only physical solutions to~(\ref{eq:intro.sh23}) --- see Refs.~\cite{ChapmanKozyreff2009,KozyreffChapman2006}. The profiles along these two branches are even-symmetric under spatial reflection (i.e., they are invariant under the reversibility transformation). The branch associated with $\phi=0$ includes profiles with a local maximum in $u$ at the midpoint, and the branch associated with $\phi=\pi$ includes profiles with a local minimum in $u$ at the midpoint. We refer to the former branch as $\brLzero$ and the latter as $\brLone$. These two branches of even-symmetric localised states persist to finite amplitude where they undergo homoclinic snaking --- a sequence of saddle-node bifurcations that cause the branches to intertwine as they oscillate back and forth across a parameter range called the snaking or pinning region~\cite{BurkeKnobloch2006,WoodsChampneys1999}. The resulting bifurcation structure, shown in figure~\ref{fig:sh23}, also includes a sequence of so-called rung branches which emerge from $\brLzero$ and $\brLone$ in pitchfork bifurcations located near the saddle-node bifurcations and cross-link the two snaking branches. The profiles that make up $\brLzero$ and $\brLone$ resemble several wavelengths of a uniform amplitude pattern connected to the trivial state by a symmetrically related pair of fronts. The profiles on the rung branches are similar to those on the snaking branches, but they are not symmetric. These profiles spontaneously break the reversibility symmetry of~(\ref{eq:intro.sh23}), so each point on the rung branches in figure~\ref{fig:sh23} actually represents two different profiles with identical norm which are related to each other by reflection in $x$.

The linear stability of the various stationary profiles is determined by linearizing~(\ref{eq:intro.sh23}) about the state. The even-symmetric localised states from $\brLzero$ and $\brLone$ are unstable at onset and change stability at each saddle-node bifurcation; profiles from the segments of the snaking branches that slant `up and to the right' on the bifurcation digram in figure~\ref{fig:sh23} are stable and those that slant `up and to the left' are unstable. All of the asymmetric profiles from the rungs are unstable.

The variational and conservative properties of the Swift-Hohenberg equation help considerably in understanding these localised states and the associated snaking bifurcation structure~\cite{BurkeKnobloch2007}. For example, the fact that the spatial dynamics associated with~(\ref{eq:intro.sh23}) is conservative determines the wavelength (i.e., the spatial period) of the pattern within the localised states. At fixed $r$, there typically exists an entire family of stationary, spatially periodic patterns $u_\mathrm{P}(x;k)$ parameterised by the wavenumber $k$. The particular pattern that is selected to appear within the localised state must lie in the level set $H=0$. Figure~\ref{fig:sh23} includes the branch $\brP$ of spatially periodic states defined by $H=0$, and the pattern wavenumber $k$ varies with $r$ along this branch to satisfy the $H=0$ constraint. Careful measurement of the numerically computed localised states confirms that this branch of patterns correctly predicts the wavenumber variation $k(r)$ within the localised states, at least when the localised states are sufficiently wide. The variational property of~(\ref{eq:intro.sh23}) is also useful in understanding the localised states. The free energy $\Lyap$ of the uniform amplitude patterns varies along $\brP$. The so-called Maxwell point M is the $r$ value at which the pattern on $\brP$ has the same free energy as the trivial state --- i.e., it is defined by the two constraints $H=\Lyap=0$. Away from the Maxwell point, there is a free energy mismatch between the pattern within the localised state and the trivial background state, with the pattern energetically favoured above the Maxwell point and the trivial flat state favoured below. This free energy mismatch is countered by the energy associated with the fronts, which pins the fronts to the interior pattern. The stationary localised states within the snaking region correspond to critical points of the free energy landscape (local minima for the stable localised states, local maxima and saddles for the unstable localised states from the snaking branches and the rungs). Outside the snaking region the free energy mismatch is sufficiently large to eliminate the critical points, de-pinning the fronts and forcing them to drift. The temporal dynamics of localised initial conditions at values of $r$ outside the snaking region support this description. The Maxwell point therefore serves as an organizing center for the snaking structure and the localised states.

It is clear that the special properties of the Swift-Hohenberg equation confer additional properties which, while extremely useful in understanding the bifurcation structure of localised states, may not be valid generically. The generic situation in which Turing instability occurs in one dimension is spatially left-right symmetric and hence the spatial dynamical system is reversible, but there is no additional requirement either for the spatial dynamical system to be conservative or for the temporal dynamics to be variational. In fact, many systems without these properties exhibit localised states which undergo homoclinic snaking, such as the complex Ginzburg-Landau equation with either 2:1 or 1:1 resonant forcing~\cite{MaETAL2010}. Similar behaviour is also observed in systems in higher spatial dimension, such as binary fluid convection in 2D~\cite{BatisteETAL2006} and plane Couette flow in 3D~\cite{SchneiderETAL2010}. Therefore it is clearly of substantial interest to understand the dynamics in cases where reversibility is preserved but the variational/conservative structure is lost. In this paper our specific motivation is from work by Kozyreff and Tlidi~\cite{KozyreffTlidi2007} who argue that in a particular double limit of small subcriticality and small critical wavenumber for the pattern forming instability, an extended version of the Swift-Hohenberg equation is necessary to capture the generic pattern-forming behaviour:
\begin{equation}\label{eq:intro.extended}
 \partial_t u = r u - \left( 1 + \partial_{xx}\right)^2 u + b u^2 - u^3 + \alpha (\partial_x u)^2 + \beta u \partial_{xx} u \, .
\end{equation}
This equation is an extension of~(\ref{eq:intro.sh23}) in which, due to the double limit, quadratic nonlinearities containing two spatial derivatives can be established to be of the same asymptotic order as the other usual terms. Like~(\ref{eq:intro.sh23}), equation~(\ref{eq:intro.extended}) is reversible; it is also equivariant under the following inversion involving the parameters: $(x,u; r, b, \alpha, \beta) \rightarrow (x,-u; r, -b, -\alpha, -\beta)$. Equation~(\ref{eq:intro.extended}) remains variational and spatially conservative for $\alpha = \beta / 2$, but in the generic case $\alpha \neq \beta /2$ it loses these properties. 
On a technical level, we note that~(\ref{eq:intro.extended}) is a semilinear PDE; extending the right-hand side further to include generic nonlinear terms containing four spatial derivatives would yield a PDE that was quasilinear but no longer semilinear.

The dynamics and structure of localised states in~(\ref{eq:intro.extended}) are formidably complicated. In this paper, we provide an initial investigation that focuses on a few key issues. Firstly, in \S\ref{sec:nf} we present the normal form analysis of the extended Swift-Hohenberg equation~(\ref{eq:intro.extended}) near the pattern forming instability. A lengthy multiple-scales calculation allows us to compute the normal form coefficients explicitly in terms of the parameters in~(\ref{eq:intro.extended}), which in turn enables us to establish the different regimes for the dynamics of this equation. Sections~\ref{sec:snaking} and~\ref{sec:onepeak} contain numerical investigations of large amplitude localised states far from onset. In \S\ref{sec:snaking} we show that localised states in the extended Swift-Hohenberg equation~(\ref{eq:intro.extended}) exhibit homoclinic snaking, though the details are not exactly the same as in the usual Swift-Hohenberg equation~(\ref{eq:intro.sh23}). In \S\ref{sec:onepeak} we examine the existence and stability of one-peak localised states. We find that there is both travelling wave (`drift') instability and standing wave instability of the single-peak state, and we discuss how these new bifurcations arise. Section~\ref{sec:conclusions} concludes.

\section{Normal form coefficients}\label{sec:nf}

\begin{figure}
\begin{center}
 \resizebox{4in}{!}{\includegraphics{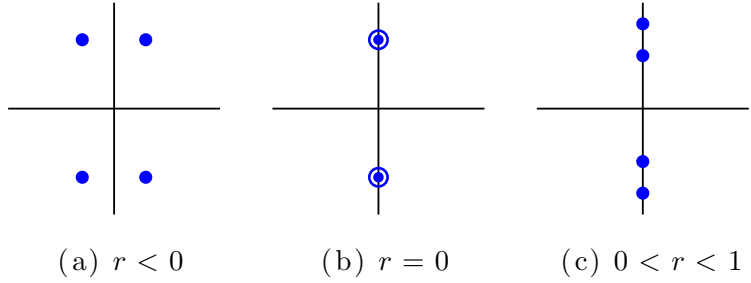}}
\end{center}
\caption{Spatial eigenvalues of the trivial state $u=0$ of the Swift-Hohenberg equation~(\ref{eq:intro.sh23}), plotted in the complex plane. (a) For $r<0$ the eigenvalues form two complex conjugate non-real pairs with real parts of equal magnitudes and opposite signs. (b) For $r=0$ the eigenvalues are purely imaginary and of multiplicity two. (c) For $0<r<1$ the eigenvalues form two purely imaginary pairs. The same eigenvalue structure applies to the normal form of the Hamiltonian-Hopf bifurcation~(\ref{eq:nf.nf}) in: (a) $\mu<0$, (b) $\mu=0$, and (c) $\mu>0$.}
\label{fig:HHeigs}
\end{figure}

In this section we use normal form analysis to examine the initial pattern forming instability of the extended Swift-Hohenberg equation~(\ref{eq:intro.extended}) at $r=0$. The time-independent version of this equation forms a fourth-order spatial dynamical system. The four spatial eigenvalues of the fixed point associated with the trivial state $u=0$ are given by $\{ \pm (\sqrt{r}-1)^{1/2}, \pm (-\sqrt{r}-1)^{1/2} \}$. In $r<0$ the four spatial eigenvalues form a complex quartet and the origin is hyperbolic. At $r=0$ the spatial eigenvalues collide pairwise on the imaginary axis, and in $0<r<1$ they form two imaginary pairs so that the origin is a center. Figure~\ref{fig:HHeigs} shows the behaviour in the spatial eigenvalues near $r=0$, which is characteristic of the Hamiltonian-Hopf bifurcation. The dynamics in the normal form of this bifurcation are well understood, and the goal of this section is to classify the dynamics of the extended Swift-Hohenberg equation~(\ref{eq:intro.extended}) by deriving the relationship between the coefficients in the normal form and the parameters in~(\ref{eq:intro.extended}). A straightforward but lengthy method to determine this relationship is to explicitly transform the spatial dynamical system associated with~(\ref{eq:intro.extended}) into normal form. This involves introducing a four-dimensional coordinate system which reproduces the spatial dynamics of~(\ref{eq:intro.extended}), then performing an appropriate linear transformation followed by a sequence of nonlinear near-identity transformations to match the normal form order by order. We choose an alternative method, based on Ref.~\cite{Grimshaw1994} and later expanded on in Ref.~\cite{BurkeKnoblochDCDS}, which involves reducing~(\ref{eq:intro.extended}) to an amplitude equation and comparing this to a suitable scaled reduction of the normal form. The values of the normal form coefficients in terms of the parameters from~(\ref{eq:intro.extended}) can be read off by comparing the two reduced equations.

We begin in \S\ref{sec:nf.nf} with a summary of the normal form for the Hamiltonian-Hopf bifurcation, paying particular attention to heteroclinic and homoclinic orbits because these correspond to fronts and localised states in~(\ref{eq:intro.extended}). In \S\ref{sec:nf.nfscale} we describe the scalings that reduce the normal form~(\ref{eq:nf.nf}) to a Ginzburg-Landau equation. In \S\ref{sec:nf.shscale} we describe a similar reduction for~(\ref{eq:intro.extended}), and use this to derive the relationship between the normal form coefficients and the parameters in~(\ref{eq:intro.extended}). In \S\ref{sec:nf.geometry} we describe geometrical features of this relationship.

\subsection{The Hamiltonian-Hopf normal form} \label{sec:nf.nf}

The normal form for the Hamiltonian-Hopf bifurcation is~\cite{IoossPeroueme1993,WoodsChampneys1999}:
\begin{subequations} \label{eq:nf.nf}
\begin{align}
 A' & = i A + B + i A P(\mu; y,w) \, , \\
 B' & = i B + i B P(\mu; y,w) + A Q(\mu; y,w) \, ,
\end{align}
\end{subequations}
where $A(x)$ and $B(x)$ are complex variables, $y \equiv |A|^2$ and $w \equiv \frac{i}{2} (A \bar{B} - \bar{A} B)$, the overbar refers to complex conjugation, and prime denotes differentiation with respect to $x$ in the context of our spatial dynamical system. The parameter $\mu$ is an unfolding parameter analogous to $r$ in~(\ref{eq:intro.extended}). The normal form is equivariant with respect to the reversibility transformation $(x,A,B) \rightarrow (-x,\bar{A}, -\bar{B})$. The terms $P(\mu;y,w)$ and $Q(\mu;y,w)$ are polynomials with real coefficients. We consider only the first few terms in $P$ and $Q$, linear in $\mu$ and up to quadratic order in $y$ and $w$:
\begin{subequations} \label{eq:nf.PQ}
\begin{align}
 P(\mu; y,w) &= p_1 \mu + p_2 y + p_3 w + p_4 y^2 + p_5 w y + p_6 w^2 \, , \\
 Q(\mu; y,w) &= -q_1 \mu + q_2 y + q_3 w + q_4 y^2 + q_5 w y + q_6 w^2 \, .
\end{align}
\end{subequations}
Without loss of generality we can choose the sign of $\mu$ to ensure $q_1>0$. With this choice, the fixed point $A=B=0$ is hyperbolic in $\mu<0$ and a center in $\mu>0$, so the bifurcation is oriented so that the spatial eigenvalues follow those shown in figure~\ref{fig:HHeigs}. The nonlinear behaviour of the unfolding near $\mu=0$ depends crucially on two of the normal form coefficients, $q_2$ and $q_4$. The sign of $q_2$ determines the criticality of the bifurcation of small amplitude periodic orbits at $\mu=0$: this bifurcation is supercritical for $q_2>0$ and subcritical for $q_2<0$. The role of $q_4$ is more subtle, as discussed below. We remark that if $q_2$ is $\order(1)$ then the role of $q_4$ is formally irrelevant to describing the dynamics near $\mu=0$. However, the calculation that follows is performed in the neighborhood of the codimension-2 point $(\mu,q_2)=(0,0)$ and so the sign of $q_4$ becomes highly relevant.

Space does not permit a complete discussion and analysis of the normal form dynamics. We present a very brief summary that is sufficient to highlight the differences in terms of possible orbits homoclinic to zero. We refer the interested reader to the analysis contained in Refs.~\cite{DiasIooss1996,Iooss1997,IoossPeroueme1993,WoodsChampneys1999} for further details.

The analysis of the dynamics of the normal form is made considerable easier by the existence of two conserved quantities: $w$ (defined above) and $h \equiv |B|^2 - \int_0^{y} Q(\mu; s, w) ds$. Within the level set $w=h=0$, which includes the fixed point $A=B=0$, the dynamics reduces to a second-order nonlinear oscillator which can be conveniently written in the form:
\begin{equation}
 \left( \frac{dy}{dx} \right)^2 + f(y) = 0 \, ,
\end{equation}
where
\begin{equation}
  f(y) = 4 q_1 \mu y^2 - 2 q_2 y^3 - \frac{4}{3} q_4 y^4 \, ,
\end{equation}
and by definition only $y \geq 0$ contains physically relevant solution trajectories.

\begin{figure}
\begin{center}
 \resizebox{6in}{!}{\includegraphics{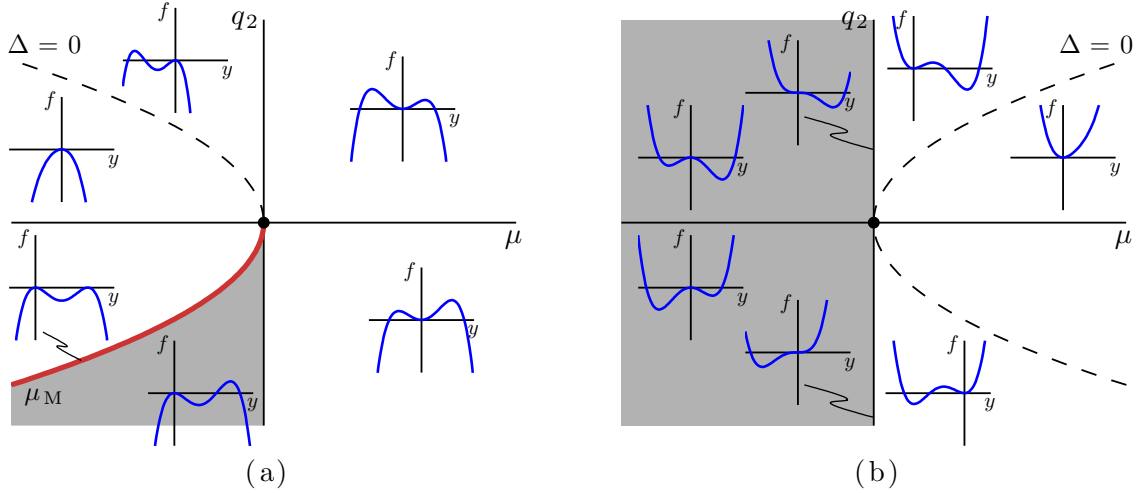}}
\end{center}
\caption{Schematic summary of the regimes of normal form dynamics in the $(\mu,q_2)$ plane, in the cases (a) $q_4>0$ and (b) $q_4<0$, after figures 2 and 3 of Ref.~\cite{WoodsChampneys1999}. Orbits homoclinic to the origin exist in the shaded regions. The dashed curves indicate $\Delta=0$, where the quantity $\Delta \equiv q_2^2/4 + 4 q_1 q_4 \mu / 3$ is the discriminant of $f(y)/y^2$. The solid line in (a) labelled $\mu_\mathrm{M}$ indicates the Maxwell point.}
\label{fig:q2q4}
\end{figure}

Figure~\ref{fig:q2q4} indicates the shape of the potential function $f(y)$ for the different sign combinations of $\mu$, $q_2$ and $q_4$.  In the case $q_4 > 0$, figure~\ref{fig:q2q4}(a) shows that orbits homoclinic to the origin exist only in $q_2<0$ and only for $\mu_\mathrm{M} \leq \mu \leq 0$, where $\mu_\mathrm{M}$ is defined by the condition that the discriminant of $f(y)/y^2$ vanishes. At $\mu_\mathrm{M}$ the potential $f(y)$ has a doubly-degenerate zero at a nontrivial value of $y$ corresponding to a periodic orbit within the $w=h=0$ level set; there is also a heteroclinic orbit (i.e., a spatial front) connecting the origin to this periodic state. The codimension-one point $\mu_\mathrm{M}$ is analogous to the Maxwell point in~(\ref{eq:intro.sh23}). While homoclinic orbits are certainly found throughout the shaded region in figure~\ref{fig:q2q4}(a), the multiplicity of localised states associated with homoclinic snaking is associated with the heteroclinic orbits along $\mu_\mathrm{M}$.

Turning to the case $q_4<0$, illustrated in figure~\ref{fig:q2q4}(b), we see that homoclinic orbits exist throughout $\mu<0$ for $q_2$ of either sign. The homoclinic orbits for $q_2<0$ are born at small amplitude as $\mu$ decreases through zero, whereas the homoclinic orbits that exist near $\mu=0$ when $q_2>0$ exist at finite amplitude. Moreover, exactly at $\mu=0$ these homoclinic orbits have algebraically decaying tails since at $\mu=0$ the origin is non-hyperbolic. In addition, Dias and Iooss~\cite{DiasIooss1996} have shown that, for $q_4<0$ and in the regime where $\mu>0$ and $q_2>0$, there exist orbits that are homoclinic to the periodic solutions.

\subsection{Scaling the normal form} \label{sec:nf.nfscale}

We wish to reduce the normal form~(\ref{eq:nf.nf}) to an amplitude equation valid near the bifurcation at $\mu=0$. The calculations that follow are simplified by focussing on the regime of mild criticality, hence we restrict ourselves to the behaviour near the codimension-two point $(\mu,q_2)=(0,0)$.

We begin by introducing a small parameter $\epsilon \ll 1$ and rescaling the parameters as $\mu=\epsilon^4 \hat{\mu}$ and $q_2 = \epsilon^2 \hat{q}_2$. Next, we define the large spatial scale $X=\epsilon^2 x$ and write $(A,B) = ( \epsilon \tilde{A}(X), \epsilon^3 \tilde{B}(X)) e^{i x}$. Dropping the tildes, the polynomials $P$ and $Q$ in~(\ref{eq:nf.PQ}) become
\begin{subequations}
\begin{align}
 P &= \epsilon^2 p_2 |A|^2 + \order(\epsilon^4)  \, , \\
 Q &= -\epsilon^4 q_1 \hat{\mu} + \epsilon^4 \hat{q}_2 |A|^2 + \epsilon^4 q_3 \frac{i}{2} (A \bar{B} - \bar{A} B) + \epsilon^4 q_4 |A|^4 + \order(\epsilon^6) \, ,
\end{align}
\end{subequations}
and so the normal form~(\ref{eq:nf.nf}) becomes
\begin{subequations}
\begin{align}
 \epsilon^3 A' &= \epsilon^3 B + i \epsilon A \left[ \epsilon^2 p_2 |A|^2 \right] + \order(\epsilon^5) \, , \label{eq:nf.nfep_A} \\
 \epsilon^5 B' &= i \epsilon^3 B \left[ \epsilon^2 p_2 |A|^2 \right] + \epsilon A \left[ -\epsilon^4 q_1 \hat{\mu} + \epsilon^4 \hat{q}_2 |A|^2 + \epsilon^4 q_3 \frac{i}{2} \left( A \bar{B} - \bar{A} B \right) + \epsilon^4 q_4 |A|^4 \right] + \order(\epsilon^7) \, . \label{eq:nf.nfep_B}
\end{align}
\end{subequations}
Rearranging~(\ref{eq:nf.nfep_A}) enables us to write $B$ in terms of $A$:
\begin{equation} \label{eq:nf.B}
 B = A' - i p_2 A |A|^2 + \order(\epsilon^2) \, .
\end{equation}
After differentiating with respect to $X$, this gives:
\begin{equation} \label{eq:nf.Bprime1}
 B' = A'' - i p_2 \left( 2 A' |A|^2 + A^2 \bar{A}' \right) + \order (\epsilon^2) \, .
\end{equation}
Note that~(\ref{eq:nf.nfep_B}) gives an alternate (and independent) asymptotic expression for $B'$:
\begin{align}
 B' &= i p_2 B |A|^2 - q_1 \hat{\mu} A + \hat{q}_2 A|A|^2 + q_3 \frac{i}{2} A (A \bar{B} - \bar{A} B) + q_4 A |A|^4 + \order (\epsilon^2) \notag \\
\label{eq:nf.Bprime2}  &= i p_2 \big( A' - i p_2 A |A|^2 \big) |A|^2
- q_1 \hat{\mu} A + \hat{q}_2 A|A|^2 \\ \notag & \quad + q_3 \frac{i}{2} A^2 (\bar{A}'
  + i p_2 \bar{A} |A|^2) - q_3 \frac{i}{2} |A|^2 (A' - i p_2 A |A|^2) + q_4 A |A|^4 + \order (\epsilon^2) \, ,
\end{align}
after eliminating factors of $B$ using~(\ref{eq:nf.B}). Equating (\ref{eq:nf.Bprime1}) and (\ref{eq:nf.Bprime2}) leads to a time-independent Ginzburg-Landau-type equation for $A(X)$:
\begin{align}
 A'' &= -q_1 \hat{\mu} A + \hat{q}_2 A |A|^2 + i \Big( p_2+\frac{1}{2}
  q_3 \Big) A^2 \bar{A}' + i \Big( 3
  p_2-\frac{1}{2} q_3 \Big) A' |A|^2 \label{eq:nf.GLnf} \\
  & \qquad\notag + (q_4 - q_3 p_2 + p_2^2 ) A|A|^4 + \order(\epsilon^2) \, .
\end{align}
This equation is an approximation of the normal form~(\ref{eq:nf.nf}) valid in the neighborhood of the codimension-2 point $(\mu,q_2)=(0,0)$.

\subsection{Scaling of the Swift-Hohenberg equation, and matching} \label{sec:nf.shscale}

We now return to the extended Swift-Hohenberg equation~(\ref{eq:intro.extended}), and introduce a scaling that leads to an analogous version of~(\ref{eq:nf.GLnf}) valid near the pattern forming instability at $r=0$. Again, we focus on a neighborhood of the codimension-2 point where the pattern forming bifurcation changes criticality. We anticipate that the condition of mild criticality will demand that a particular combination of the quadratic coefficients $(b, \alpha, \beta)$ is small, and we allow this combination to emerge naturally in the calculation.

The reduction procedure for~(\ref{eq:intro.extended}) involves taking appropriate scalings for the parameters and expanding $u(x,t)$ as a sum of Fourier modes multiplied by amplitudes that depend on long spatial and temporal scales. We introduce a small parameter $\epsilon \ll 1$ and rescale the parameters as:

\begin{equation} \label{eq:nf.scaleSHparams}
 r = \epsilon^4 \hat{\mu} \, , \quad b = b_0 + \epsilon^2 \hat{b} \, , \quad \alpha = \alpha_0 + \epsilon^2 \hat{\alpha} \, , \quad \beta = \beta_0 + \epsilon^2 \hat{\beta} \, ,
\end{equation}
where $b_0$, $\alpha_0$, $\beta_0$ are the values (to be determined) of the quadratic coefficients that correspond to $q_2=0$, and $\hat{\mu}$, $\hat{b}$, $\hat{\alpha}$, $\hat{\beta}$ are all $\order(1)$. Next, we define the large spatial scale $X = \epsilon^2 x$ and long time scale $T=\epsilon^4 t$, and propose the following ansatz for solutions to~(\ref{eq:intro.extended}):
\begin{equation} \label{eq:nf.SHansatz}
 u(x,t) = \epsilon^2 \Theta + \left[ \epsilon A e^{i x} + \epsilon^2 B e^{2 i x}  + \epsilon^3 C e^{3 i x} + \epsilon^4 D e^{4 i x} + \cc \right] + \order\left( \epsilon^4 \right) \, ,
\end{equation}
where the amplitudes $\Theta$, $A$, $B$, $C$, $D$ are functions of $X$ and $T$ and are all $\order(1)$, the higher order terms in $\epsilon$ take the form $\epsilon^n e^{n i x} + \cc$ for $n\geq5$, and $\cc$ denotes complex conjugation of the terms preceeding it within the brackets. Note that the usual approach is to proceed order by order in $\epsilon$, bringing in appropriate modes at each order to avoid secular terms. We choose instead to include in~(\ref{eq:nf.SHansatz}) all the Fourier modes necessary, with scalings motivated in obvious ways from the number of quadratic interactions necessary to produce a term with the appropriate Fourier dependence. We proceed by substituting~(\ref{eq:nf.scaleSHparams}) and~(\ref{eq:nf.SHansatz}) into~(\ref{eq:intro.extended}) and collecting terms with the same Fourier dependence $e^{n i x}$, keeping careful track throughout this procedure of the size of the largest neglected terms.
The results for $n=0,1,2,3$ are:
\begin{subequations}
{\allowdisplaybreaks
\begin{align}
 n=0: \label{eq:nf.SHn0} && 0 & = -\epsilon^2 \Theta + b_0 \big( 2 \epsilon^2 |A|^2
 + \epsilon^4 \Theta^2 + 2 \epsilon^4 |B|^2 \big) \\ &&& \quad - 3
 \big( 2 \epsilon^4 \Theta |A|^2 + \epsilon^4 \bar{A}^2 B + \epsilon^4
 A^2 \bar{B} \big) \notag \\ &&& \quad + \alpha_0 \big( 2
 \epsilon^2 |A|^2 + 8 \epsilon^4 |B|^2 + 2 \epsilon^2 (i \epsilon^2 A
 \partial_X \bar{A} - i \epsilon^2 (\partial_X A) \bar{A}) \big) \notag \\ &&& \quad +
 \beta_0 \big( -2\epsilon^2 |A|^2 - 8 \epsilon^4 |B|^2 + 2
 \epsilon^2 ( i \epsilon^2 (\partial_X A) \bar{A} - i \epsilon^2 A \partial_X \bar{A}) \big)
 \notag \\ &&& \quad + 2 \epsilon^4  \hat{b} |A|^2
  + 2 \epsilon^4 \hat{\alpha} |A|^2
 - 2 \epsilon^4 \hat{\beta} |A|^2
  + \order(\epsilon^6) \notag \\
 n=1: \label{eq:nf.SHn1}  && \epsilon^5 \partial_T A & = \epsilon^5 \hat{\mu} A + \epsilon^5 \partial_{XX} A + b_0
 \big( 2 \epsilon^3 \Theta A + 2 \epsilon^3 \bar{A} B + 2 \epsilon^5
 \bar{B} C \big) \\ &&& \quad - 3 \big(
 \epsilon^3 A |A|^2 + \epsilon^5 \Theta^2 A + \epsilon^5 \bar{A}^2 C +
 2 \epsilon^5 \Theta \bar{A} B + 2 \epsilon^5 A |B|^2 \big) \notag \\
 &&& \quad + \alpha_0 \big( 4 \epsilon^3 \bar{A} B + 12
 \epsilon^5 \bar{B} C + 2 i \epsilon^5 (\partial_X \Theta) A + 4 i
 \epsilon^5 (\partial_X \bar{A}) B - 2 i \epsilon^5 \bar{A} \partial_X B  \big) \notag \\ &&&
 \quad + \beta_0 \big( - \epsilon^3 \Theta A - 5 \epsilon^3
 \bar{A} B - 13 \epsilon^5 \bar{B} C + 2 i \epsilon^5
 \Theta \partial_X A + 4 i \epsilon^5 \bar{A} \partial_X B - 2 i \epsilon^5 (\partial_X \bar{A}) B)
 \notag \\ &&& \quad  + 2 \epsilon^5 \hat{b} \big( \Theta A + \bar{A} B \big) +
 4 \epsilon^5 \hat{\alpha} \bar{A} B -
 \epsilon^5 \hat{\beta} \big( \Theta A + 5 
 \bar{A} B \big) + \order(\epsilon^7) \notag \\
 n=2: \label{eq:nf.SHn2} && 0 & = -9 \epsilon^2 B + 24 i \epsilon^4 \partial_X B + b_0 \big(
 \epsilon^2 A^2 + 2 \epsilon^4 \Theta B + 2 \epsilon^4 \bar{A} C \big)
 \\ &&& \quad - 3 \big( \epsilon^4 \Theta A^2 + 2 \epsilon^4 |A|^2 B
 \big) \notag \\ &&& \quad + \alpha_0 \big( -\epsilon^2 A^2 + 6
 \epsilon^4 \bar{A} C + 2 i \epsilon^4 A \partial_X A \big) \notag
 \\ &&& \quad + \beta_0 \big( -\epsilon^2 A^2 - 4 \epsilon^4
 \Theta B - 10 \epsilon^4 \bar{A} C + 2 i \epsilon^4  A \partial_X A
 \big) \notag \notag \\ &&& \quad + \epsilon^4 \hat{b} A^2 - \epsilon^4 \hat{\alpha} 
 A^2  - 2 \epsilon^4 \hat{\beta} A^2 +
 \order(\epsilon^6) \notag \\
 n=3: \label{eq:nf.SHn3} && 0 & = -64\epsilon^3 C + b_0 \big( 2 \epsilon^3 A B \big) - \big( \epsilon^3 A^3 \big) - \alpha_0 \big( 4 \epsilon^3 A B \big) - \beta_0 \big( 5 \epsilon^3 A B \big) + \order(\epsilon^5) \, .
\end{align}
}
\end{subequations}
The $n=1$ terms in~(\ref{eq:nf.SHn1}) are the equivalent of the `solvability conditions' found using other methods, and these are clearly the terms in which we are ultimately interested. However, we need first to solve~(\ref{eq:nf.SHn0}),~(\ref{eq:nf.SHn2}) and~(\ref{eq:nf.SHn3}) to express the amplitudes $\Theta$, $B$ and $C$ in terms of the principal amplitude $A$. We do this iteratively by writing
\begin{equation} \label{eq:nf.ThetaBC}
 \Theta = \Theta_0 + \epsilon^2 \Theta_2 + \order(\epsilon^4) \, , \quad
 B = B_0 + \epsilon^2 B_2 + \order(\epsilon^4) \, , \quad
 C = C_0 + \order(\epsilon^2) \, ,
\end{equation}
and substituting these into~(\ref{eq:nf.SHn0}),~(\ref{eq:nf.SHn2}) and~(\ref{eq:nf.SHn3}). The leading order relations are:
\begin{equation} \label{eq:nf.ThetaBC_0}
 \Theta_0 = c_1 |A|^2 \, , \quad
 B_0 = c_2 A^2 \, , \quad
 C_0 = c_3 A^3 \, ,
\end{equation}
where we define the parameter combinations
\begin{subequations} \label{eq:nf.c1c2c3}
\begin{align}
 c_1 &= 2 \left( b_0 + \alpha_0 - \beta_0 \right) \, , \\
 c_2 &= \frac{1}{9} \left( b_0 - \alpha_0 - \beta_0 \right) \, , \\
 c_3 & = \frac{1}{64} \left[ \left( 2 b_0 - 4 \alpha_0 - 5\beta_0 \right)c_2 - 1\right] \, .
\end{align}
\end{subequations}
At next order, (\ref{eq:nf.SHn0}) and~(\ref{eq:nf.SHn2}) give
\begin{subequations} \label{eq:sh.ThetaBC_2}
\begin{align}
 \Theta_2 &= \hat{c}_1 |A|^2 + c_4 |A|^4 + i c_5 \left( A \partial_X \bar{A} - (\partial_X A) \bar{A}\right) \, , \\
 B_2 &= \hat{c}_2 A^2 + c_6 A^2 |A|^2 + i c_7 A \partial_X A \, ,
\end{align}
\end{subequations}
where we define the further parameter combinations
\begin{subequations}
\begin{align}
 \hat{c}_1 &= 2 \left( \hat{b} + \hat{\alpha} - \hat{\beta} \right) = \frac{\partial c_1}{\partial b_0} \hat{b} + \frac{\partial c_1}{\partial \alpha_0} \hat{\alpha} + \frac{\partial c_1}{\partial \beta_0} \hat{\beta} \, , \\
 \hat{c}_2 &= \frac{1}{9} \left( \hat{b} - \hat{\alpha} - \hat{\beta} \right) = \frac{\partial c_2}{\partial b_0} \hat{b} + \frac{\partial c_2}{\partial \alpha_0} \hat{\alpha} + \frac{\partial c_2}{\partial \beta_0} \hat{\beta} \, , \\
 c_4 & = b_0 (c_1^2 + 2 c_2^2) - 6 (c_1+c_2) + 8 c_2^2
 (\alpha_0 - \beta_0) \, , \\
 c_5 & = 2 (\alpha_0 - \beta_0) \, , \\
 c_6 & = \frac{1}{9} \left[ 2 c_1 c_2 (b_0 - 2\beta_0) + 2 c_3 (b + 3 \alpha_0 - 5 \beta_0) - 3 (c_1+2 c_2) \right] \, , \\
 c_7 & = \frac{2}{9} \left[ 24 c_2 + (\alpha_0 + \beta_0) \right]  \, ,
\end{align}
\end{subequations}
and we use the `hat' notation in~(\ref{eq:sh.ThetaBC_2}) to indicate that $\hat{c}_1$ and $\hat{c}_2$ can be interpreted as higher order contributions to coefficients of terms that already appeared at leading order in~(\ref{eq:nf.ThetaBC_0}).

We now turn to consideration of~(\ref{eq:nf.SHn1}). As it stands, this equation contains both $\order(\epsilon^3)$ and $\order(\epsilon^5)$ terms. The problem of mixed asymptotic orders disappears if we demand that the $\order(\epsilon^3)$ terms vanish, and this defines the desired relation between the quadratic coefficients of~(\ref{eq:intro.extended}) at the codimension-2 point. Inserting~(\ref{eq:nf.ThetaBC}) and~(\ref{eq:nf.ThetaBC_0}) into~(\ref{eq:nf.SHn1}), we see that the $\order(\epsilon^3)$ terms are:
\begin{align}
 0 & = b_0 \left( 2 \Theta_0 A + 2 \bar{A} B_0 \right) - 3 A |A|^2 + 4 \alpha_0 \bar{A} B_0 - \beta_0 \left( \Theta_0 A + 5 \bar{A} B_0 \right) \, , \notag \\
  & = \left[ 3 - c_1 \left( 2 b_0 - \beta_0 \right) - c_2 \left( 2 b_0 + 4 \alpha_0 - 5 \beta_0 \right) \right] A|A|^2 \, .
\end{align}
The bracketed quantity determines the criticality of the bifurcation at $r=0$. Based on the preceding discussion of the normal form, we set
\begin{equation} \label{eq:nf.q2}
 q_2(b, \alpha, \beta) = \frac{1}{4} \left[ 3 - 2 \left( b + \alpha - \beta \right) \left( 2 b - \beta \right) - \frac{1}{9} \left( b - \alpha - \beta \right) \left( 2 b + 4 \alpha - 5 \beta \right) \right] \, ,
\end{equation}
so that at the codimension-2 point we have $q_2(b_0,\alpha_0,\beta_0)=0$ and the $\order(\epsilon^3)$ terms in~(\ref{eq:nf.SHn1}) vanish. The remaining terms in~(\ref{eq:nf.SHn1}) are $\order(\epsilon^5)$ and, after substituting and tidying up, these give a differential equation for the principal amplitude $A(X,T)$:
\begin{equation} \label{eq:nf.GLsh}
 \partial_T A = \hat{\mu} A + 4 \partial_{XX} A - 4 \hat{q}_2 A |A|^2 + i c_8 |A|^2 \partial_X A + i c_9 A^2 \partial_X \bar{A} + c_{10} A |A|^4 \, ,
\end{equation}
where:
\begin{subequations}
\begin{align}
 \hat{q}_2 &= \left. \left( \frac{\partial q_2}{\partial b} \hat{b} + \frac{\partial q_2}{\partial \alpha} \hat{\alpha} + \frac{\partial q_2}{\partial \beta} \hat{\beta} \right) \right|_{(b, \alpha, \beta) = (b_0, \alpha_0, \beta_0)} \, , \\
 c_8 &= b_0 (-2 c_5 + 2 c_7) + \alpha_0 (4 c_7 + 2 c_1 - 4 c_2) + \beta_0 (c_5 - 5 c_7 + 2 c_1 + 8 c_2) \, , \\
 c_9 &= b_0 (2 c_5) + \alpha_0 (2 c_1 + 4 c_2) + \beta_0 (-c_5 - 2 c_2) \, , \\
 c_{10} &= b_0 (2 c_4 + 2 c_6 + 2 c_2 c_3) - 3 (c_1^2 + c_3 + 2 c_1 c_2 + 2 c_2^2)  \notag \\ & \qquad {} + \alpha_0 (4 c_6 + 12 c_2 c_3) + \beta_0 (-c_4 - 5 c_6 - 13 c_2 c_3) \, .
\end{align}
\end{subequations}
Equation~(\ref{eq:nf.GLsh}) is the Ginzburg-Landau approximation to the extended Swift-Hohenberg equation~(\ref{eq:intro.extended}) valid near onset in the regimes of small criticality --- i.e., near the codimension-two point \mbox{$(r,q_2(b,\alpha, \beta)) = (0,0)$}.

Both~(\ref{eq:nf.GLnf}) and the time-independent version of~(\ref{eq:nf.GLsh}) describe the spatial evolution of the amplitude $A$, so we can compare these two equations term by term in order to identify the normal form coefficients in terms of the parameters of~(\ref{eq:intro.extended}). This procedure gives the following:
\begin{subequations}
\begin{align}
    q_1  &=-1/4 \, , \\
    p_2  &=-\frac{1}{16} (c_8 + c_9) \, , \\ 
    q_3  &= \frac{1}{8} (c_8 - 3 c_9)  \, , \\
    q_4  &= \frac{1}{256} \big(-3 c_8^2 + 2 c_8 c_9 + 5 c_9^2 \big) -
    \frac{1}{4} c_{10} \label{eq:nf.q4} \, .
\end{align}
\end{subequations}
For our purposes, it is sufficient to calculate only these normal form coefficients. The remaining normal form coefficients, as well as those associated with higher order terms that were omitted from~(\ref{eq:nf.PQ}), may be obtained by carrying out these expansions to higher order in $\epsilon$.

In the limit $\alpha = \beta = 0$, the above expressions reduce to the familiar normal form results for the usual Swift-Hohenberg equation~(\ref{eq:intro.sh23}). In this limit,~(\ref{eq:nf.q2}) becomes $q_2(b_1,0,0) = (27-38 b^2)/36$. The condition $q_2=0$ gives $b^2 = 27/38$, and at this value of $b$ equation~(\ref{eq:nf.q4}) gives $q_4 |_{q_2=0} = 2202/361$.

\subsection{Geometry of the parameter dependence of $q_2$ and $q_4$} \label{sec:nf.geometry}

In this section we describe briefly some geometrical features of the dependence of the normal form coefficients $q_2$ and $q_4$ on the parameters $b$, $\alpha$ and $\beta$. 

\begin{figure}
\begin{center}
 \resizebox{6in}{!}{\includegraphics{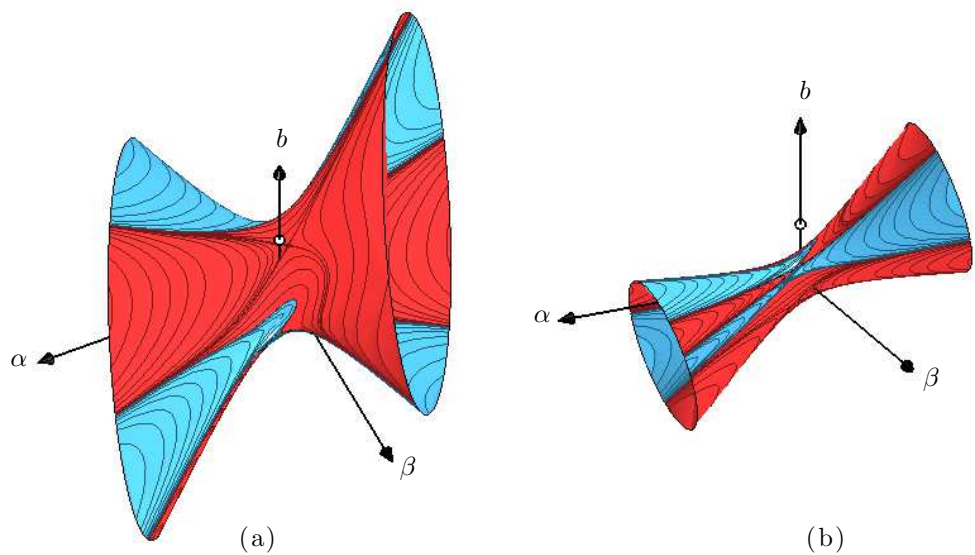}}
\end{center}
\caption{Two views of the surface $q_2=0$ plotted in the $\mathbb{R}^3$ parameter space of the quadratic coefficients $(b, \alpha, \beta)$, with contours showing level sets of $q_4$. The surface is coloured according the value of $q_4$: red and blue respectively correspond to positive and negative values of $q_4$. The open circle ($\circ$) marks the point $(b, \alpha, \beta) = (1.8,0,0)$. The surface can be viewed in the accompanying movie file (\texttt{http://math.bu.edu/people/jb/Research/extendedSH/Figure{\char`\_}4{\char`\_}movie.mpg}).}
\label{fig:hyperboloid}
\end{figure}

Expression~(\ref{eq:nf.q2}) for $q_2$ is a homogeneous quadratic function of $b_1$, $\alpha$ and $\beta$ apart from a constant term. After some simplification (\ref{eq:nf.q2}) can be written in terms of a quadratic form as:
\begin{equation} \label{eq:q2M}
 q_2 \left( b, \alpha, \beta \right) = \frac{1}{36} \left( 27 - \mathbf{p}^T M \mathbf{p} \right) \, ,
\end{equation}
where $\mathbf{p} = [b \ \ \alpha \ \ \beta]^{\mathrm{T}}$ denotes the vector of coefficients of the quadratic terms in~(\ref{eq:intro.extended}) and
\begin{equation} \label{eq:nf.M}
 M = \left[ \begin{matrix} 38 & 17 & -61/2 \\ 17 & -4 & -17/2 \\ -61/2 & -17/2 & 23 \end{matrix} \right]
\end{equation}
is a $3 \times 3$ symmetric matrix with eigenvalues $\lambda_1 \simeq 66.82$, $\lambda_2 \simeq 0.5113$, $\lambda_3 \simeq -10.33$. The condition $q_2=0$ gives
\begin{equation} \label{eq:nf.hyperboloid}
  27 = \mathbf{p}^T M \mathbf{p} \, ,
\end{equation}
which can be interpreted as an equation for a two-dimensional surface in the $\mathbb{R}^3$ parameter space of the quadratic coefficients $(b, \alpha, \beta)$. The signs of the eigenvalues of $M$ indicate that this surface is a `hyperboloid of 1-sheet'. Figure~\ref{fig:hyperboloid} shows two views of the surface~(\ref{eq:nf.hyperboloid}). The normal form coefficient $q_2$ is positive inside the hyperboloid (i.e., in the connected component of $\mathbb{R}^3$ that contains the origin) and negative outside. The parameter inversion symmetry of~(\ref{eq:intro.extended}) is apparent in the reflection symmetry of the surface.

The quantity $q_4$ is important only when $q_2$ is small, so it is natural to consider the sign of $q_4$ on the surface $q_2=0$. In figure~\ref{fig:hyperboloid} the regions in which $q_4$, as given by~(\ref{eq:nf.q4}), is positive or negative are coloured red and blue, respectively. These are labelled directly in figure~\ref{fig:hyperboloid}c which also shows, indicated by the red and white dot, the location of the usual Swift-Hohenberg equation~(\ref{eq:intro.sh23}) for which $\alpha = \beta = 0$. This is confirmation of the result mentioned above that $q_4 |_{q_2=0}>0$ in~(\ref{eq:intro.sh23}). Moreover, it is clear from figure~\ref{fig:hyperboloid} that there are substantial regions of parameter space in which $q_4<0$ and so there is ample motivation to study the dynamics in both cases of the sign of $q_4$. 

\begin{figure}
\begin{center}
 \resizebox{6in}{!}{\includegraphics{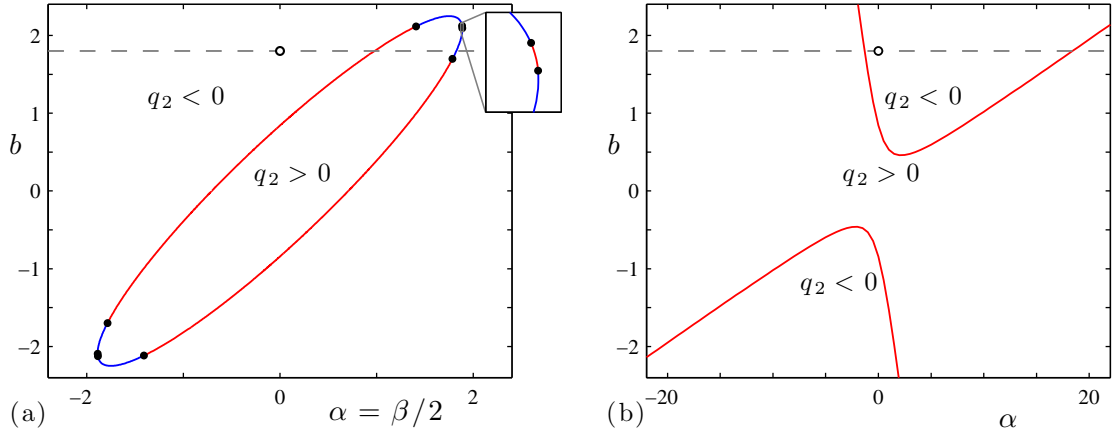}}
\end{center}
\caption{Two planar sections of the $q_2=0$ surface in figure~\ref{fig:hyperboloid}: (a) $\alpha=\beta/2$, where equation~(\ref{eq:intro.extended}) is variational in $t$ and conservative in $x$, and (b) $\beta=0$, which is the focus of \S\ref{sec:snaking} and \S\ref{sec:onepeak}. The colouring of the $q_2=0$ curve indicates the sign of $q_4$, with $q_4>0$ red and $q_4<0$ blue. Solid dots ($\bullet$) mark the points with $q_4=0$. The dashed line in each frame indicates $b=1.8$, and the open circle ($\circ$) marks the point $(b, \alpha, \beta) = (1.8,0,0)$. The inset in (a) shows a small $q_4>0$ range along the curve.}
\label{fig:slice}
\end{figure}


Recall that the extended Swift-Hohenberg equation~(\ref{eq:intro.extended}) is variational in $t$ and conservative in $x$ when $\alpha = \beta/2$. This corresponds to a planar section of the three-dimensional parameter space defined above. Figure~\ref{fig:slice}(a) shows the corresponding section of the $q_2=0$ surface from figure~\ref{fig:hyperboloid} using coordinates $(\alpha, b)$. Figure~\ref{fig:slice}(b) shows a second planar section of the $q_2=0$ surface from figure~\ref{fig:hyperboloid}, at $\beta=0$. The behaviour of~(\ref{eq:intro.extended}) at $\beta=0$ is considered in more detail in the following section.

\section{Homoclinic snaking in the extended Swift-Hohenberg equation} \label{sec:snaking}

In this section we examine the behaviour of localised states in the extended Swift-Hohenberg equation~(\ref{eq:intro.extended}). The primary goal is simply to establish (numerically) that the homoclinic snaking behaviour present in~(\ref{eq:intro.sh23}) persists with the inclusion of the non-variational and non-conservative terms. We also point out several ways in which the localised states in~(\ref{eq:intro.extended}) differ from those in~(\ref{eq:intro.sh23}). For simplicity, we fix $b=1.8$ throughout this section, and use the behaviour of~(\ref{eq:intro.sh23}) at this value of $b$ (as shown in figure~\ref{fig:sh23}) as a point of reference. We make use of~(\ref{eq:nf.q2}) to focus on values of the quadratic parameters $(b,\alpha, \beta)$ for which $q_2<0$ and $|q_2|$ is $\order(1)$ --- i.e., the highly subcritical regime, where homoclinic snaking is prominent --- and we use the software package AUTO~\cite{AUTO} to trace out branches of localised states in the primary bifurcation parameter $r$.

\begin{figure}
\begin{center}
 \resizebox{5in}{!}{\includegraphics{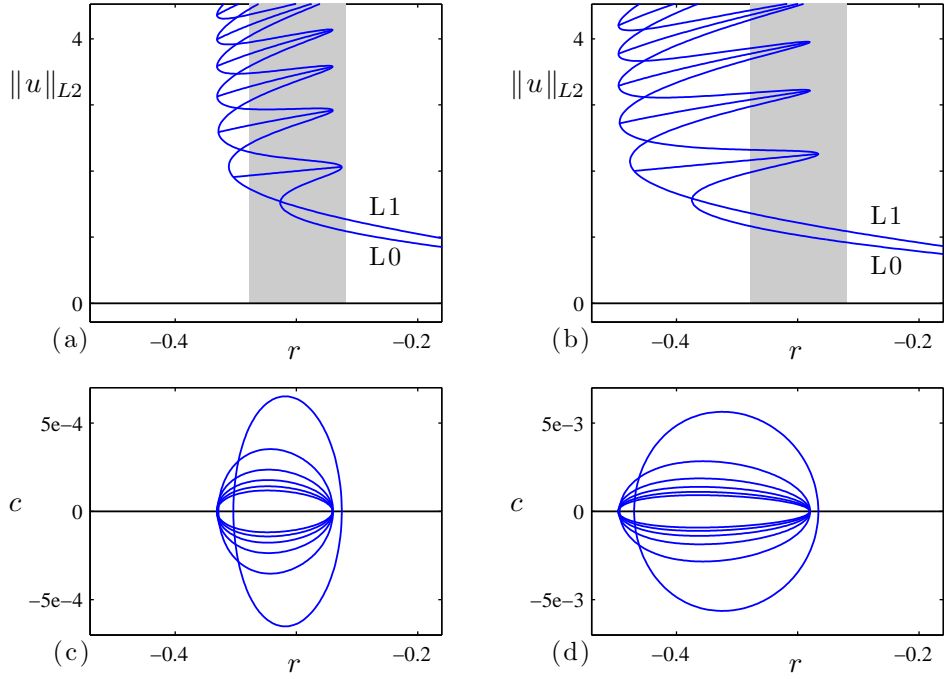}}
\end{center}
\caption{Bifurcation diagrams of localised states at (a) $\alpha=0.1$ and (b) $\alpha=0.5$, with other parameters fixed at $b=1.8$ and $\beta=0$. The two snaking branches $\brLzero$ and $\brLone$ contain even-symmetric, steady, spatially localised profiles. The profiles on the rung branches are asymmetric and drift in $x$ at constant velocities $\pm c$. Temporal stability of solutions is not indicated. The shading indicates the snaking region at $\alpha=0$ (see figure~\ref{fig:sh23}). The lower frames show the drift velocity of the asymmetric localised states from the lowest six rung branches.}
\label{fig:snakeA1}
\end{figure}

To begin, we examine the effect of the new terms in~(\ref{eq:intro.extended}) by increasing $\alpha$ from $\alpha=0$ while leaving $\beta=0$ fixed. Frames (a) and (b) of figures~\ref{fig:snakeA1} show the bifurcation diagrams of localised states at $\alpha=0.1$ and $\alpha=0.5$, respectively. The behaviour of stationary, even-symmetric localised states is qualitatively the same as that shown in figure~\ref{fig:sh23} at $\alpha = 0$. These localised states are organised in a pair of intertwined snaking branches, which we continue to label $\brLzero$ and $\brLone$. The saddle-node bifurcations on the snaking branches line up asymptotically (except perhaps the lowest two on the $\brLzero$ branch) to two $r$ values which define the snaking region. The shaded region in each frame indicates the snaking region at $\alpha=0$, so increasing $\alpha$ increases the width of the snaking region and also shifts it to more negative values of $r$. Note however that in the case of equation~(\ref{eq:intro.extended}) we can not define a Maxwell point within this snaking region. Moreover, unlike in figure~\ref{fig:sh23}, we do not indicate the stability of any localised states in figure~\ref{fig:snakeA1}. A complete description of the stability of these states is beyond the scope of this article, but we do present a limited discussion of stability in \S\ref{sec:onepeak} where we focus on one-peak solutions between the first and second saddle-node bifurcations on the $\brLzero$ branch.


At the fixed values of $b$ and $\beta$ used in figure~\ref{fig:snakeA1}, equation~(\ref{eq:nf.q2}) for $q_2$ has a zero at $\alpha \simeq -1.306$, with $q_2<0$ for values of $\alpha$ larger than this --- see figure~\ref{fig:slice}(b). From~(\ref{eq:nf.q4}) we find that $q_4>0$ at this point. The normal form therefore predicts that the behaviour in the neighborhood of the codimension-two point $(r,\alpha) = (0,-1.306)$ follows figure~{\ref{fig:q2q4}}(a), just like the familiar case of equation~(\ref{eq:intro.sh23}) in the neighborhood of $(r,b) = (0,\sqrt{27/38})$. For values of $\alpha$ slightly greater than $\alpha \simeq -1.306$, the snaking region is  exponentially narrow and close to $r=0$. At larger $\alpha$, the quantity $q_2$ decreases to an $\order(1)$ negative value, and the system becomes more subcritical. This causes the width of the snaking region to increase and shift to more negative $r$ values. This is consistent with the snaking behaviour shown at the various $\alpha$ values in figure~\ref{fig:snakeA1}.

We note that, at the fixed values $b_1=1.8$ and $\beta=0$ used in figure~\ref{fig:snakeA1}, equation~(\ref{eq:nf.q2}) for $q_2$ also has a second zero at $\alpha \simeq 18.4$ --- see figure~\ref{fig:slice}(b). The pattern forming instability at $r=0$ is subcritical in $\alpha \in [-1.306,18.4]$. Moreover, at this second codimension-two point $(r,\alpha) \simeq (0,18.4)$ we again have $q_4>0$ and so the normal form analysis suggests that nearby the behaviour of the system is similar to that shown in figure~\ref{fig:q2q4}(a). Understanding the behaviour of~(\ref{eq:intro.extended}) in the neighborhood of this second codimension-two point, and its influence on, and relationship to, the localised states associated with the first codimension-two point at $(r, \alpha) \simeq (0, -1.306)$, remains an open problem.

The bifurcation diagrams in frames (a) and (b) of figure~\ref{fig:snakeA1} also include the rung branches which cross-link the two snaking branches. The asymmetric localised states on the rungs drift in $x$ since this is generic in non-variational systems. The rung branches shown in figure~\ref{fig:snakeA1} consist of profiles of fixed shape that drift at constant velocity --- i.e., they are travelling wave solutions of the form $u(x,t) = u(x-c t)$. The drift velocity $c$ for the various rungs is plotted in the lower frames in figure~\ref{fig:snakeA1}. The drift velocity varies with $r$ along each rung, but approaches $c=0$ at the endpoints that mark the secondary pitchfork bifurcations where the rungs connect to the snaking branches. Recall that each point on a rung actually includes two different profiles related by reflection symmetry. These profiles drift with equal speed in opposite directions, so figure~\ref{fig:snakeA1}(c) includes two segments (one in $c>0$ and one in $c<0$) for each rung in figure~\ref{fig:snakeA1}(a), and likewise for figures~\ref{fig:snakeA1}(d) and~\ref{fig:snakeA1}(b). There are two trends in $c$ that are apparent in figure~\ref{fig:snakeA1}. First, at fixed $\alpha$ the drift velocity varies from one rung to the next: the narrow localised states on the lower rungs tend to travel faster than the wider localised states on the upper rungs. Second, by comparing frames (c) and (d) it appears that the drift speed tends to increase with $\alpha$, at least for small $\alpha$. Intuitively this is to be expected as one moves away from the variational case $\alpha=0$. Similar trends have been observed in drift speeds of asymmetric localised profiles in other systems~\cite{GibsonETAL}.

\begin{figure}
\begin{center}
 \resizebox{4in}{!}{\includegraphics{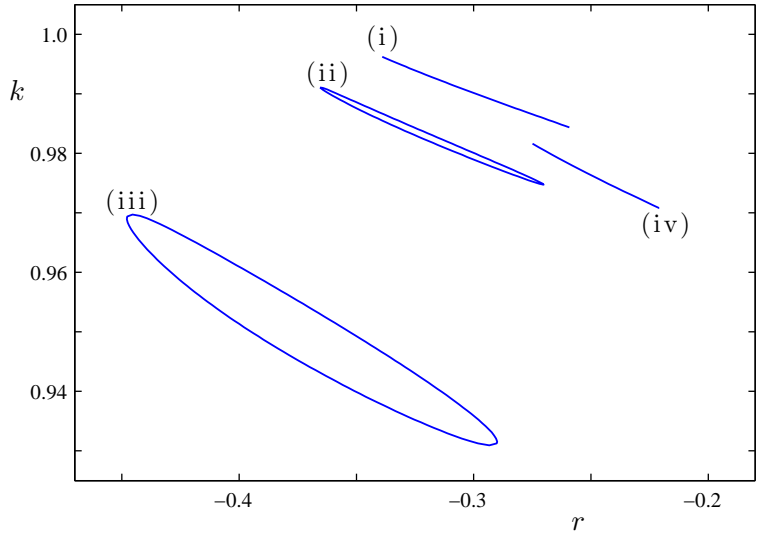}}
\end{center}
\caption{Wavenumber $k$ of the pattern within the localised states, measured along the snaking branches for several values of the quadratic coefficients: (i) $(b, \alpha, \beta) = (1.8,0,0)$, as in figure~\ref{fig:sh23}; (ii) $(b, \alpha, \beta) = (1.8,0.1,0)$, as in figure~\ref{fig:snakeA1}(a); (iii) $(b, \alpha, \beta) = (1.8,0.5,0)$, as in figure~\ref{fig:snakeA1}(b); and (iv) $(b, \alpha, \beta) = (1.8,0.1,0.2)$. For the parameter values in (i) and (iv), the system is conservative in $x$ because $\alpha=\beta/2$, and the wavenumber selection $k(r)$ is a single-valued function of $r$. For the parameter values in (ii) and (iv), the system is non-conservative in $x$ and the wavenumber selection $k(r)$ is multivalued, tracing out an isola. The wavenumbers shown in the figure only apply to the wide localised states taken from far up the respective snaking branches.}
\label{fig:wavenumber}
\end{figure}

The new terms in~(\ref{eq:intro.extended}) also cause a qualitative change in the wavenumber selection relative to that found in~(\ref{eq:intro.sh23}). Figure~\ref{fig:wavenumber} shows the measured wavenumber $k$ of the pattern within the localised states as a function of $r$ for one complete back-and-forth cycle across the pinning region, for several $(\alpha,\beta)$ pairs at fixed $b=1.8$. In each case, the measurement involves wide profiles from far up the snaking branches where the wavenumber is well defined. Though the measurements are made using one particular back-and-forth pair of segments from one particular snaking branch, they are characteristic of all the wide localised states far up both the $\brLzero$ and $\brLone$ snaking branches. For the narrow localised states lower on the snaking branches, the wavenumber is poorly defined.

The curve labelled (i) in figure~(\ref{fig:wavenumber}) shows the wavenumber variation of the localised states along the snaking branches from figure~\ref{fig:sh23}, at $\alpha=\beta=0$. Recall that in this case the spatial dynamics is conservative so the patterns within the localised states must satisfy $H=0$. This constraint $H=0$ defines the branch $\brP$ of patterns included in figure~\ref{fig:sh23}, and determines the variation in wavenumber $k$ along this branch. Thus $k(r)$ for the patterns is defined uniquely over a range of $r$ which includes the snaking region of localised states; the measured wavenumber on the curve (i) in figure~\ref{fig:wavenumber} is the segment of this broader $k(r)$ curve that lies within the snaking region. In particular, we note that the wavenumber variation of the stable localised states along the segments of the snaking branches that slant `up and to the right' is identical to the wavenumber variation of the unstable localised states along the segments that slant `up and to the left'. Two profiles at the same $r$ value on consecutive segments include identical patterns in their interior, and differ only in the shape of the fronts that connect the pattern to the flat background.  Two distinct fronts are created in a saddle-node bifurcation at the left edge of the pinning region, and merge in a second saddle-node bifurcation at the right edge~\cite{BeckETAL2009}. The $H=0$ constraint forces the fronts to approach the same spatially periodic orbit.

In the absence of a spatially conserved quantity $H$, the saddle-node bifurcation of fronts creates a pair of fronts that approach two different periodic orbits. The wavenumber variation along the segments of the snaking branches that slant `up and to the right' is therefore generically different from that along the segments that slant `up and to the left', though of course they match at the saddle-node bifurcations at the edge of the snaking region. The curves (ii) and (iii) in figure~\ref{fig:wavenumber} show the wavenumber measured along the snaking branches from frames (a) and (b) of figure~\ref{fig:snakeA1}. In these examples, the spatial dynamics is non-conservative in $x$ and the wavenumber along the snaking branches traces out a loop. The larger wavenumber corresponds to segments of the snaking branches that slant `up and to the right', and the lower wavenumber to segments that slant `up and to the left', at least for these two examples. We note that the splitting of $k(r)$ into a loop has also been observed in other spatially non-conservative systems, such as plane Couette flow~\cite{GibsonETAL}.

The curve labelled (iv) in figure~\ref{fig:wavenumber} shows the wavenumber variation at $\alpha=0.1$ and $\beta=0.2$. As $\alpha=\beta/2$, this corresponds to a spatially conservative limit of~(\ref{eq:intro.extended}). As such, the wavenumber variation is again a single-valued function of $r$. The splitting of the wavenumber selected by the localised states is therefore a measure of how spatially non-conservative is the system, though the high-precision numerical measurements required to observe the splitting may make this impractical as a diagnostic tool.

\section{One peak solutions} \label{sec:onepeak}

In this section we investigate another aspect of how the non-variational and non-conservative terms effect the behaviour of ~(\ref{eq:intro.extended}). We focus on the stability of one-peak localised states since this enables us to show clearly the emergence of two new kinds of temporal instabilities: drift and standing oscillations. The one-peak localised states occur on the segment of the $\brLzero$ branch between the first and second saddle-node bifurcations. It is convenient to introduce the notation $\segonepeak$ to refer to this segment. We organise the results by fixing $b=1.8$ and $\beta=0$, and examining the states on $\segonepeak$ at increasing values of $\alpha$.

Throughout this section, we report on the linear stability of the one-peak states. As usual, the linear stability of a stationary solution $u_0(x)$ to~(\ref{eq:intro.extended}) is determined by substituting $u(x,t) = u_0(x) + \epsilon U(x) e^{\sigma t}$ into~(\ref{eq:intro.extended}) and ignoring $\order(\epsilon^2)$ terms. The mode $U(x)$ and growth rate $\sigma$ satisfy the linear eigenvalue equation $\sigma U = \mathcal{L}[u_0, \partial_x] U$ where the linearised operator associated with the right side of~(\ref{eq:intro.extended}) is:
\begin{equation}
 \mathcal{L}[u_0,\partial_x] = r - \left( 1+\partial_{xx} \right)^2 + 2 b u_0 - 3 u_0^2 + 2 \alpha (\partial_x u_0) \partial_x + \beta \left( (\partial_{xx} u_0) + u_0 \partial_{xx} \right) \, .
\end{equation}
Spatial translation invariance of~(\ref{eq:intro.extended}) implies that any steady solution always has a neutrally stable ($\sigma=0$) Goldstone mode $U_\mathrm{G}(x) = \partial_x u_0(x)$. The one-peak solutions of interest here are even-symmetric, and the associated Goldstone mode is odd-symmetric. In what follows, we report the existence of unstable ($\mathrm{Re}\,\sigma > 0$) modes using the notation $[m_r,m_c,n_r,n_c]$, where $m_r$ ($m_c$) is the number of even-symmetric modes with real (complex) eigenvalues $\sigma$, and $n_r$ ($n_c$) is the number of odd-symmetric modes with real (complex) eigenvalues. Stable solutions are labelled $[0,0,0,0]$. Note that both $m_c$ and $n_c$ must be even; in what follows, we never observe unstable odd-symmetric modes with complex eigenvalues so $n_c=0$ throughout.

\begin{figure}
\begin{center}
 \resizebox{5in}{!}{\includegraphics{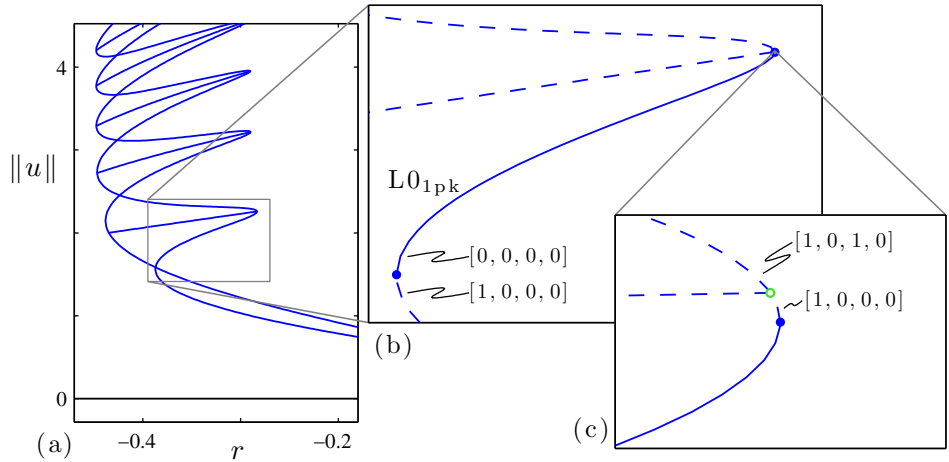}}
\end{center}
\caption{Bifurcation diagram of one-peak states for $(b,\alpha, \beta) = (1.8,0.5,0)$. (a) Snaking diagram of localised states, with no indication of stability. (b) Detail of the $\segonepeak$ segment; solid (dashed) curves indicate stable (unstable) states. The number and symmetry of unstable modes is also indicated, using the notation defined in the text. (c) Detail of the right edge of $\segonepeak$. Solid dots ($\bullet$) indicate saddle-node bifurcations; the open circle ($\circ$) indicates the secondary bifurcation to the rung branch. One-peak solutions are stable over the entire $\segonepeak$ segment.}
\label{fig:onepeakA05}
\end{figure}

\begin{figure}
\begin{center}
 \resizebox{6in}{!}{\includegraphics{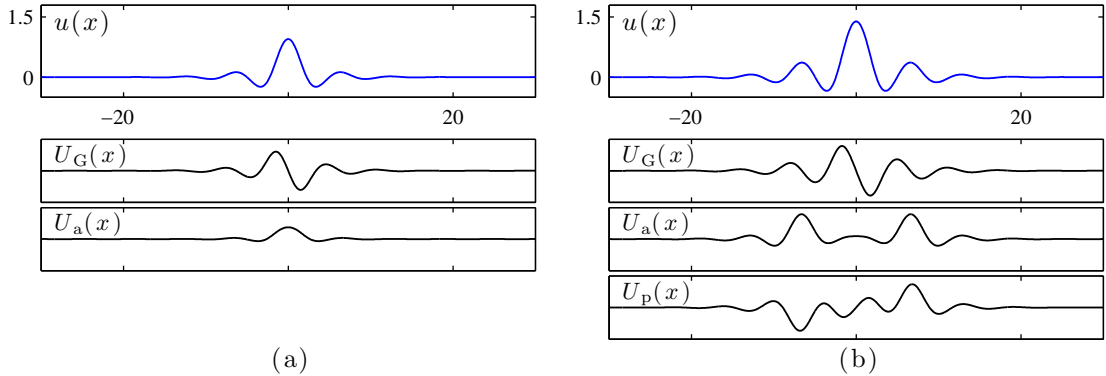}}
\end{center}
\caption{One-peak profiles and eigenfunctions from figure~\ref{fig:onepeakA05} at $(b,\alpha,\beta) = (1.8,0.5,0)$, (a) near the left edge of $\segonepeak$ and (b) near the right edge of $\segonepeak$. Each column includes the stationary one-peak state $u(x)$, the neutrally stable Goldstone mode $U_\mathrm{G}(x)$, and the even-symmetric amplitude mode $U_\mathrm{a}(x)$ that changes stability at the corresponding saddle-node bifurcation. The right column also includes the odd-symmetric phase mode $U_\mathrm{p}(x)$ which is associated with the lowest rung branch.}
\label{fig:onepeakA05prof}
\end{figure}

A natural starting point is the case $\alpha=0$, as the behaviour of the usual Swift-Hohenberg equation~(\ref{eq:intro.sh23}) is well-known. As shown in figure~\ref{fig:sh23}, the one-peak localised states on $\segonepeak$ are all stable. Increasing $\alpha$ initially maintains qualitatively very similar behaviour. Figure~\ref{fig:onepeakA05} shows the behaviour at $\alpha=0.5$. Figure~\ref{fig:onepeakA05}(a) shows the usual snaking behaviour of $\brLzero$ and $\brLone$. Figure~\ref{fig:onepeakA05}(b) shows a detail of the $\segonepeak$ segment, and includes stability assignments. Below the saddle-node bifurcation at the left edge of $\segonepeak$, a single even-symmetric `amplitude' mode becomes unstable. Likewise, at the saddle-node bifurcation at the right edge of $\segonepeak$ an even-symmetric amplitude mode becomes unstable. This is followed by an odd-symmetric `phase' mode which becomes unstable slightly above this saddle-node bifurcation, at the secondary bifurcation that gives rise to the lowest rung branch. Figure~\ref{fig:onepeakA05prof} shows the profiles $u(x)$ at two points along $\segonepeak$, one near the left edge and one near the right edge. This figure also includes plots of the modes that are most important in determining stability. We note in passing that, unlike in the case of the usual Swift-Hohenberg equation~(\ref{eq:intro.sh23}), the amplitude mode which passes through zero growth rate at the left edge of $\segonepeak$ is different to the one that passes through zero growth rate at the right edge. The former merges with another even-symmetric mode along and becomes complex with $\mathrm{Re}\,\sigma<0$, corresponding to a stable oscillatory mode. 

\begin{figure}
\begin{center}
 \resizebox{5in}{!}{\includegraphics{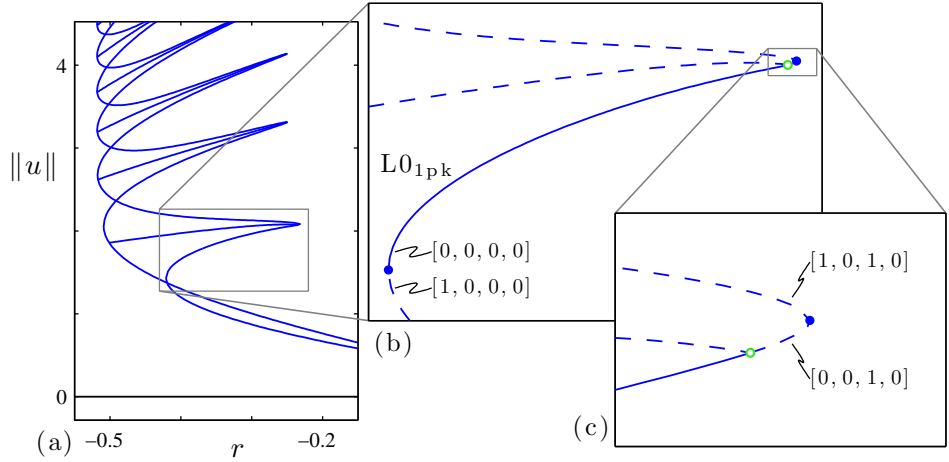}}
\end{center}
\caption{Bifurcation diagram of one-peak states for $(b,\alpha, \beta) = (1.8,1.6,0)$. Notation follows figure~\ref{fig:onepeakA05}. Stationary one-peak solutions are stable from the left edge of $\segonepeak$ up to $r \simeq -0.2358$  where the odd-symmetric mode associated with the rung loses stability.}
\label{fig:onepeakA16}
\end{figure}

\begin{figure}
\begin{center}
 \resizebox{6in}{!}{\includegraphics{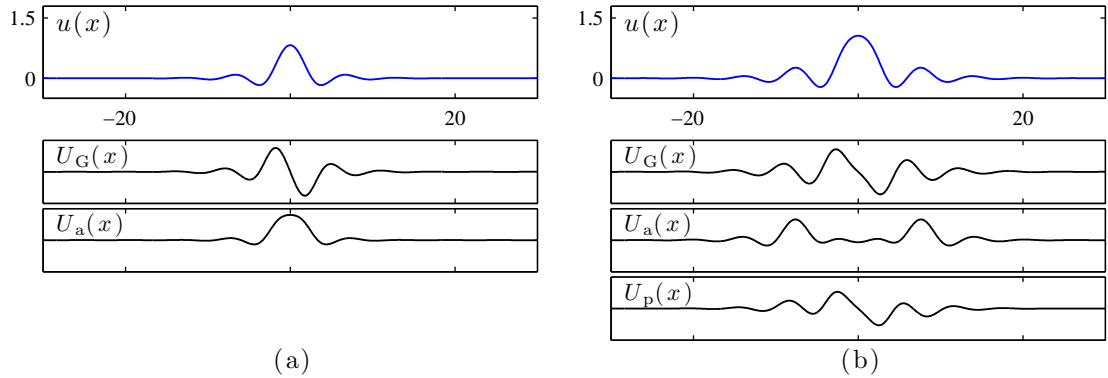}}
\end{center}
\caption{One-peak profiles and eigenfunctions from figure~\ref{fig:onepeakA16} at $(b,\alpha,\beta) = (1.8,1.6,0)$, (a) near the left edge of $\segonepeak$ and (b) near the right edge of $\segonepeak$. Notation follows figure~\ref{fig:onepeakA05prof}.}
\label{fig:onepeakA16prof}
\end{figure}

Qualitatively new behaviour occurs at $\alpha=1.6$, as shown in figures~\ref{fig:onepeakA16} and~\ref{fig:onepeakA16prof}. The bifurcation behaviour in the neighborhood of the left edge of $\segonepeak$ remains unchanged --- below this saddle-node bifurcation there is a single unstable amplitude mode which stabilises at the saddle-node bifurcation. However, the profiles on $\segonepeak$ are not stable all the way up to the saddle-node bifurcation at the right edge. Instead, a phase mode loses stability slightly before the right edge. This is the same instability as that which was previously associated with the secondary bifurcation to a rung branch, so now the rung connects to the $\brLzero$ branch slightly below the right saddle-node bifurcation. Though not labelled in the figure, we find that the profiles on the rung branch remain unstable.

\begin{figure}
\begin{center}
 \resizebox{5in}{!}{\includegraphics{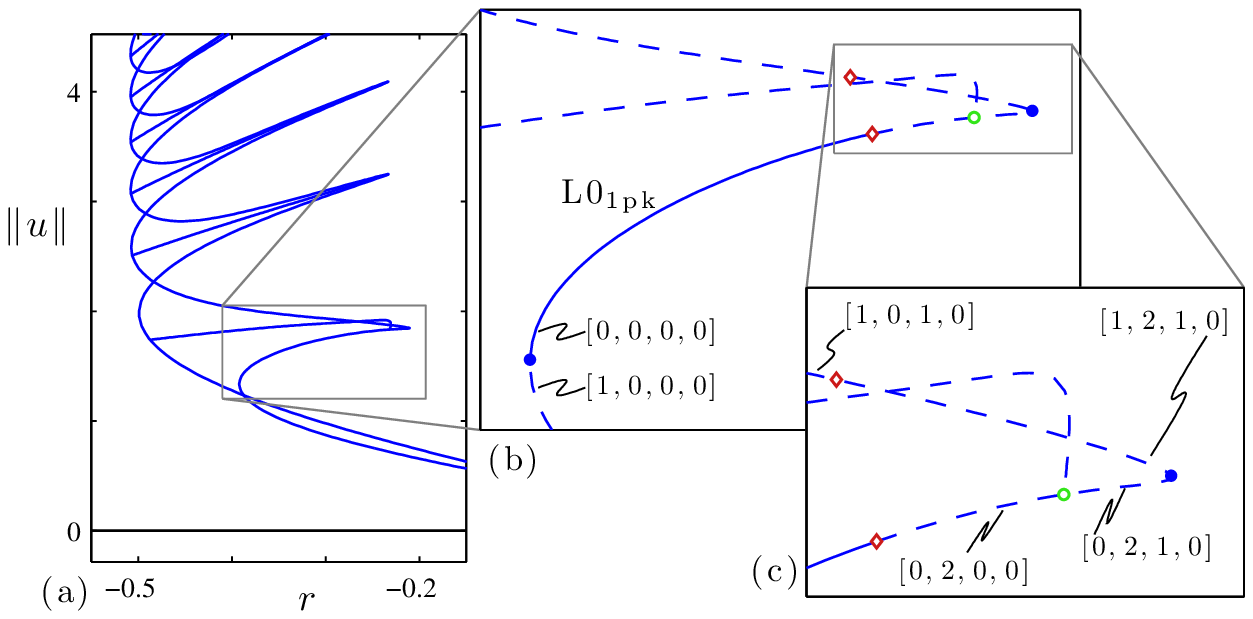}}
\end{center}
\caption{Bifurcation diagram of one-peak states for $(b,\alpha, \beta) = (1.8,2.0,0)$. Notation follows figure~\ref{fig:onepeakA05}. Stationary one-peak solutions are stable from the left edge of $\segonepeak$ up to $r \simeq -0.2682$ where an even-symmetric oscillatory mode loses stability; the oscillatory instability is marked with a diamond symbol ($\diamond$).}
\label{fig:onepeakA20}
\end{figure}

\begin{figure}
\begin{center}
 \resizebox{6in}{!}{\includegraphics{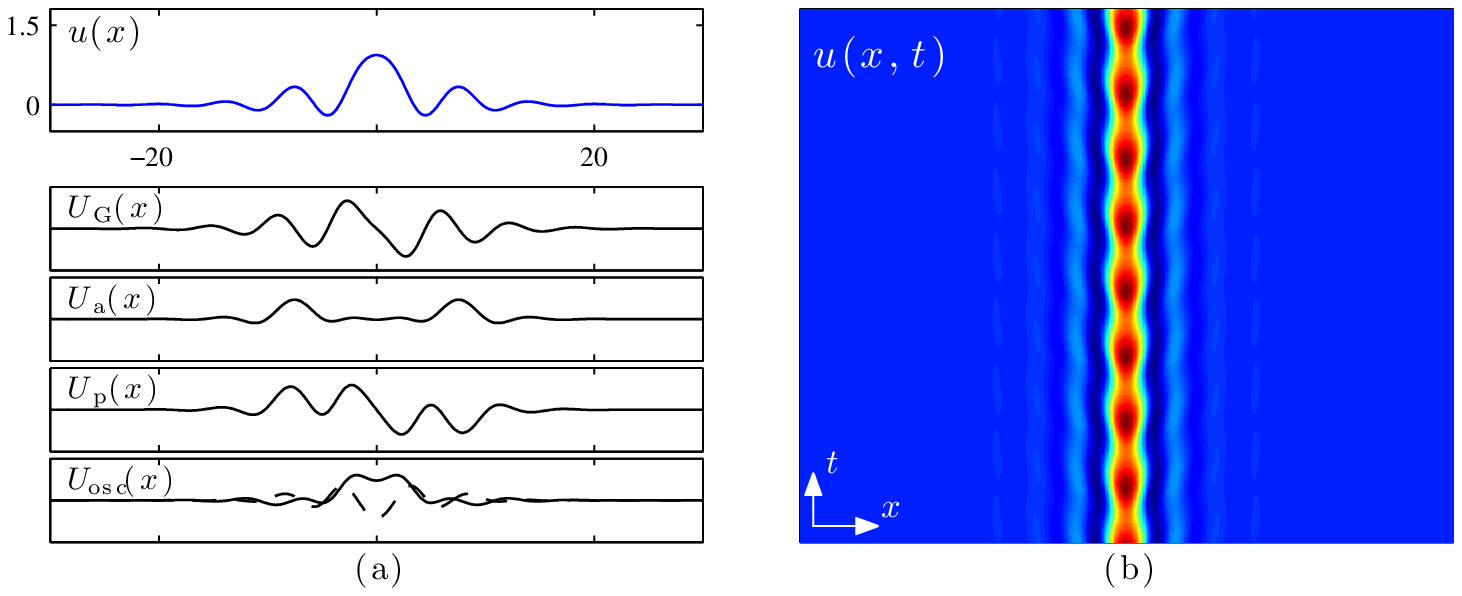}}
\end{center}
\caption{(a) One-peak profile and eigenfunctions from figure~\ref{fig:onepeakA20} at $(b,\alpha,\beta) = (1.8,2.0,0)$, at a value of $r$ slightly above the oscillatory instability on $\segonepeak$. Notation follows figure~\ref{fig:onepeakA05prof}. The eigenfunction $U_\mathrm{osc}(x)$ of the mode associated with the oscillatory instability is also included; solid and dased lines indicate the real and imaginary parts of $U_\mathrm{osc}(x)$ respectively. (b) Space-time plot of an oscillon at $(b,\alpha, \beta) = (1.8,2.0,0)$, shown on the domain $(x,t) \in [-50,50] \times [0,100]$. The $r$ value for this solution is $r = -0.2681$, which is slightly above the oscillatory instability of the one-peak solutions from figure~\ref{fig:onepeakA20}. The oscillon can be viewed in the accompanying movie file (\texttt{http://math.bu.edu/people/jb/Research/extendedSH/Figure{\char`\_}13{\char`\_}movie.mpg}).}
\label{fig:onepeakA20prof}
\end{figure}

Increasing $\alpha$ causes the secondary bifurcation to the rung branch to move further from the right edge of $\segonepeak$, but also leads to new instability. The bifurcation diagram at $\alpha = 2.0$ is shown in figure~\ref{fig:onepeakA20}. At this value of $\alpha$, the first mode to become unstable along $\segonepeak$ is an even-symmetric oscillatory mode, which loses stability at $r \simeq -0.2682$. The relevant mode $U_\mathrm{osc}(x)$ is therefore complex, and its real and imaginary parts are shown in figure~\ref{fig:onepeakA20prof}(a). This mode eventually re-stabilises much further along $\brLzero$, above the saddle-node bifurcation at the right edge of $\segonepeak$. We find that the imaginary part of the complex eigenvalue at onset is $\mathrm{Im}\,\sigma \simeq 0.55$ and that it remains nearly constant along the snaking branch. 

The stationary localised solutions on $\segonepeak$ at $r$ values above the complex instability are unstable, but may evolve to time-dependent states which remain spatially localised. Such solutions are often referred to as oscillons. Figure~\ref{fig:onepeakA20prof}(b) shows the space-time plot of one such solution at $r=-0.2681$, slightly above the initial instability on $\segonepeak$ at $r \simeq -0.2682$. The criticality of this instability is, however, difficult to confirm. Numerical results indicate that the oscillons that occur above onset are large perturbations of the unstable stationary solutions on $\segonepeak$. Furthermore, the oscillons persist as $r$ decreases below $r \simeq -0.2682$. One explanation for this might be that the oscillatory instability is subcritical, and that the branch of oscillons that emerges from this point turns around and stabilises in a saddle-node bifurcation. A second possibility is that the branch that emerges from $\segonepeak$ remains unstable, and the localised solution shown in figure~\ref{fig:onepeakA20prof}(b) lies on a separate branch of localised oscillations that emerges elsewhere in the bifurcation diagram. This second possibility occurs in the autonomous system of reaction-diffusion equations studied in Ref.~\cite{AvitabileETAL}, which includes a bifurcation structure very similar to that shown in figure~\ref{fig:onepeakA20}.

Further increase in $\alpha$ causes the oscillatory instability to invade more of the stable range of the $\segonepeak$ segment. At $\alpha \simeq 2.8$ this instability reaches the left edge of $\segonepeak$ so all stationary one-peak solutions are unstable.

\section{Conclusions} \label{sec:conclusions}

In this paper we have considered the effect of non-variational and non-conservative terms on the well-known quadratic-cubic Swift-Hohenberg equation. Such an investigation is well motivated in general by the wealth of recent work on homoclinic snaking in a variety of contexts where no variational principle exists, and in particular by the paper of Kozyreff and Tlidi~\cite{KozyreffTlidi2007} who show that~(\ref{eq:intro.extended}) arises naturally as a model equation for a long-wavelength pattern forming instability in one-dimensional dissipative systems. 

The dynamics of the extended Swift-Hohenberg equation~(\ref{eq:intro.extended}) are extremely rich and we have been content to set the scene for future work by focussing on a few aspects of the problem. Firstly we have carried out the computation of the normal form coefficients for the pattern-forming bifurcation problem. This is a key first step in the analysis of this extended equation. It allows us to use the normal form results as a guide to the kinds of localised states that should exist in different parameter regimes. One important result is the existence of substantial open regions of parameter space in which $q_4<0$.

Next we investigated numerically the effect of one of the new terms on the snaking behaviour at large amplitude, far from the initial linear instability. We found that the snaking phenomenon persists, and we pointed out aspects of snaking in the extended Swift-Hohenberg equation~(\ref{eq:intro.extended}) that differ from the well-known snaking in~(\ref{eq:intro.sh23}), including wavenumber selection within the localised states, and the behaviour of the solutions from the rungs.

Lastly we investigated the emergence of two new kinds of instability of stationary one-peak localised states that affect the parameter range over which stable one-peak states exist. These are (i) a drift instability and (ii) an oscillatory instability. The latter may lead to stable localised oscillations, but the bifurcation structure of such solutions remains unclear. We showed how the two new instabilities arise as qualitative changes to the traditional snaking bifurcation diagram. Since both instabilities are generic we expect that this analysis will provide a guide to the location and dynamics of instabilities of multi-peak localised states and quite possibly multi-pulse localised states as well. The details of these cases, as so much else, we leave for future work.

\vspace{0.2cm}
\noindent
{\bf Acknowledgement}: 
J.B. acknowledges the award of a David Crighton Fellowship from DAMTP, University of Cambridge where this work was initiated. The research of J.B. was also supported by the Center for BioDynamics at Boston University and the NSF (DMS 0602204, EMSW21-RTG). JHPD holds a University Research Fellowship from the Royal Society.

\bibliography{extendedSH}

\end{document}